\documentclass[11pt,           % in 11pt fonts
english                       % default language
]{article}

\usepackage[english,strings]{babel}     % with explicit language
\usepackage{amsmath}           % ams mathematical stuff
\usepackage[utf8]{inputenc}    % smart input of funny chars
\usepackage[T1]{fontenc}       % also for the font encoding
\usepackage{longtable}         % tables longer than one page
\usepackage{exscale}           % large summation signs in 11pt
\usepackage[final]{graphicx}   % to include pdf pictures
\usepackage[sort]{cite}        % nicer citations
\usepackage{array}             % nice tables
\usepackage{fancyhdr}          % nicer headers
\usepackage[a4paper]{geometry} % geometry of page layout
\usepackage{gitinfo2} % show info about the git
\usepackage[multiuser]{fixme}  % correction notes, warnings etc.
\usepackage{xspace}            % better spacing after macros
\usepackage{tikz}              % for commutative diagrams and stuff
\usepackage{ifdraft}           % to determine whether draft mode
\usepackage{nchairx}
\usepackage[expansion=false    % no font expansion
           ]{microtype}        % only protrusion
\usepackage[nottoc]{tocbibind} % refs and index in the toc
\usepackage[backref=page,      % backrefs in the bibliography
           final=true,         % always treat as final
           pdfpagelabels       % use pdf page labels
           ]{hyperref}         % hyperrefs are cool!

\graphicspath{{../tikz/}}

\usetikzlibrary{matrix,calc,fadings,decorations.pathreplacing,positioning,patterns}

\ifdraft{\synctex=1}{}

\fxusetheme{color}

\newcommand{\myname}{\textbf{Stefan Waldmann}}
\newcommand{\myemail}{\texttt{stefan.waldmann@mathematik.uni-wuerzburg.de}}
\newcommand{\myaddress}{Julius Maximilian University of Würzburg \\
     Department of Mathematics \\
     Chair of Mathematics X (Mathematical Physics) \\
     Emil-Fischer-Straße 31 \\
     97074 Würzburg \\
     Germany}

\author{\myname\thanks{\myemail}\\[0.5cm]\myaddress}

\newcommand{\lie}[1]{\liealg{#1}}
\newcommand{\halbnorm}[1]{\seminorm{#1}}
\newcommand{\starwick}{\mathbin{\star_{\scriptscriptstyle\mathrm{Wick}}}}
\newcommand{\starred}{\mathbin{\star_{\scriptscriptstyle\mathrm{red}}}}
\renewcommand{\basis}[1]{\mathsf{#1}}

\title{Convergence of Star Product: From Examples to a General
  Framework\thanks{Awarded with the \emph{Prix du Concour annuel 2018: On
      demande une contribution à la construction d’un cadre convergent
      pour la quantification par déformation}
    of the \emph{Académie royale des Sciences, des Lettres et des
      Beaux-Arts de Belgique.}}}

\date{January 2019}

\begin{document}

\selectlanguage{english}

%
% title page
%

\maketitle

%
% abstract
%

\begin{abstract}
    We recall some of the fundamental achievements of formal
    deformation quantization to argue that one of the most important
    remaining problems is the question of convergence. Here we discuss
    different approaches found in the literature so far. The recent
    developments of finding convergence conditions are then outlined
    in three basic examples: the Weyl star product for constant
    Poisson structures, the Gutt star product for linear Poisson
    structures, and the Wick type star product on the Poincaré disc.
\end{abstract}

\newpage

%
% Introduction: Formal Deformation Quantization
%

\section{Introduction: Formal Deformation Quantization}
\label{sec:FormalDeformationQuantization}

In the by now classical paper \cite{bayen.et.al:1978a}, Bayen,
Fr{\o}nsdal, Flato, Lichnerowicz, and Sternheimer introduced the
notion of a formal star product on a Poisson manifold as a general
method to pass from a classical mechanical system encoded by the
Poisson manifold to its corresponding quantum system. The main
difference compared to other, and more ad-hoc, quantization schemes is
the emphasis on the role of the observable algebra. The observable
algebra is constructed not as particular operators on a Hilbert space
as this is usually done. Instead, one stays with the same vector space
of smooth functions on the Poisson manifold and just changes the
commutative pointwise product into a new noncommutative product using
Planck's constant $\hbar$ as a deformation parameter. This way one
tries to model the quantum mechanical commutation relations.

While in \cite{bayen.et.al:1978a} it is shown that this leads to the
same results as expected in quantum mechanics for those systems where
alternative quantization schemes are available, the proposed scheme of
deformation quantization has several conceptual advantages: the first
is of course its \emph{vast generality} concerning the
formulation. While other quantization schemes make much more use of
specific features of the classical system, deformation quantization
can be seen as almost universal in the sense that its requirements are
virtually minimal. Only a Poisson algebra of classical observables is
needed to formulate the program of deformation quantization. Of
course, the hard part of the work consists then in actually proving
the existence (and possible classifications) of star products, but in
any case, the conceptual framework is fixed from the beginning. A
second advantage is that the \emph{physical interpretation} of the
observables is fixed from the beginning: the observables simply stay
the same elements of the same underlying vector space. Thus it is
immediately clear which quantum observable is the Hamiltonian, the
momentum etc. since they are the same as classical. It is only the
product law which changes, the correspondence of classical and quantum
observables is implemented trivially.  A third advantage of this
approach to focus on the algebra first is shared also by other
formulations of quantum theory, most notably by axiomatic quantum
field theory, see e.g. \cite{haag:1993a}: having the focus on the
algebra one can now study \emph{different representations} on Hilbert
spaces which might be needed to encode different physical situations
the system is exposed to. In quantum field theory this is a well-known
feature and difficulty, but now the same options are available in
quantum mechanical systems as well. It turns out that this is not just
a mathematical game to play but required by physical reality. As
argued in e.g.  \cite{bordemann.neumaier.pflaum.waldmann:2003a}
inequivalent, representations are needed to encode the Aharonov-Bohm
effect \cite{aharonov.bohm:1959a} in this framework.

There is of course a price to be paid in order to achieve this
generality in the formulation of the quantization problem. To
understand the difficulties we recall the precise definition of a
formal star product: let $(M, \pi)$ be a Poisson manifold, i.e. a
smooth manifold with a Poisson tensor $\pi \in \Secinfty(\Anti^2 TM)$
such that we have a Poisson bracket $\{f, g\} = \pi(\D f, \D g)$ for
smooth functions $f, g \in \Cinfty(M)$. With other words, the
commutative associative algebra $\Cinfty(M)$ of complex-valued smooth
functions becomes a Poisson algebra when equipped with this bracket
$\{\argument, \argument\}$. Then a formal star product $\star$ for
$(M, \pi)$ consists in a product law
\begin{equation}
    \label{eq:FormalStarProduct}
    \star\colon
    \Cinfty(M)[[\lambda]] \times \Cinfty(M)[[\lambda]]
    \ni (f, g)
    \; \mapsto \;
    f \star g
    \in \Cinfty(M)[[\lambda]]
\end{equation}
defined on the space of formal power series in a formal parameter
$\lambda$ with coefficients in the smooth functions such that the
following properties hold: first, $\star$ is
$\mathbb{C}[[\lambda]]$-bilinear. This already implies that there are
uniquely determined $\mathbb{C}$-bilinear maps
$C_r\colon \Cinfty(M) \times \Cinfty(M) \longrightarrow \Cinfty(M)$
such that
\begin{equation}
    \label{eq:OrdersStar}
    f \star g
    =
    \sum_{r=0}^\infty \lambda^r C_r(f, g),
\end{equation}
where the $C_r$ are extended to formal series by requiring
$\mathbb{C}[[\lambda]]$-bilinearity. The second requirement is that
$\star$ should be associative. This is in fact a truly nontrivial
feature as it results in a infinite chain of quadratic equations for
the coefficient operators $C_r$. The first orders of the associativity
condition are the associativity of $C_0$ in zeroth order, then
\begin{equation}
    \label{eq:CrAssociative}
    C_0(C_1(f, g), h) + C_1(C_0(f, g), h)
    =
    C_0(f, C_1(g, h)) + C_1(f, C_0(g, h))
\end{equation}
in first order and so on. The third requirement is the compatibility
with the semi-classical limit $\lambda = 0$. One wants
\begin{equation}
    \label{eq:SemiclassicalLimit}
    C_0(f, g) = fg
    \quad
    \textrm{and}
    \quad
    C_1(f, g) - C_1(g, f) = \I \{f, g\}
\end{equation}
for all $f, g \in \Cinfty(M)$. This way, one implements the actual
quantization condition that the product $\star$ deforms its zeroth
order $C_0$ which should be the usual product and second, the
commutator with respect to $\star$ deforms the Poisson bracket in the
next nontrivial order. Including $\I = \sqrt{-1}$ in the requirement
as in \eqref{eq:SemiclassicalLimit} then allows for the identification
\begin{equation}
    \label{eq:LambdaBecomeshbar}
    \lambda \leadsto \hbar,
\end{equation}
i.e. the formal deformation parameter $\lambda$ should correspond to
the physical Planck constant $\hbar$. The last requirements are more
of a technical nature and not too essential for the physical
interpretation. One requires $1 \star f = f = f \star 1$ for the
constant function $1$ and one requires the $C_r$ to be bidifferential
operators. It turns out that in all known constructions these two
features are automatically satisfied and will therefore not pose any
additional difficulties.

With this definition, a star product $\star$ becomes a particular
example of a formal deformation of an associative algebra in the sense
of Gerstenhaber \cite{gerstenhaber:1963a, gerstenhaber:1964a,
  gerstenhaber:1966a, gerstenhaber:1968a, gerstenhaber:1974a}. The
main difficulty of formal deformation quantization is now that the
deformation parameter $\lambda$ corresponds to the constant $\hbar$
which is one of the very fundamental constants of nature. In
particular, this constant is non-zero and \emph{not}
dimensionless. Thus it will make no sense whatsoever to speak of
smallness of $\hbar$ in an intrinsic way. In fact, one can always
choose a unit system where $\hbar$ has the numerical value $\hbar = 1$
or, say, $\hbar = 42$. Thus the convergence of the series
\eqref{eq:OrdersStar} becomes an inevitable issue which has to be
solve in order to have a physically sound model. This is the
convergence problem in deformation quantization.

Before entering the details on the convergence issue let us first
continue with some positive results in deformation quantization
illustrating that it is worth to pursue the proposed path. Indeed, on
the level of formal star products, deformation quantization is
extraordinarily successful: after several classes of examples
\cite{cahen.gutt:1981a, gutt:1983a, cahen.gutt:1983b,
  cahen.gutt:1983a, dewilde.lecomte:1983c, dewilde.lecomte:1983a} the
existence of star products on general symplectic manifolds was shown
by de Wilde and Lecomte \cite{dewilde.lecomte:1983b}. Shortly later,
Fedosov gave a different and geometrically very explicit construction
\cite{fedosov:1985a, fedosov:1986a, fedosov:1994a}, see also his book
\cite{fedosov:1996a}. Finally, Omori, Maeda and Yoshioka found yet
another construction for symplectic manifolds
\cite{omori.maeda.yoshioka:1991a}. The classification of formal star
products on symplectic manifolds was obtained be several groups in
slightly different formulations in \cite{deligne:1995a,
  gutt.rawnsley:1999a, bertelson.cahen.gutt:1997a, nest.tsygan:1995a,
  nest.tsygan:1995b, fedosov:1996a, weinstein.xu:1998a,
  neumaier:2002a}. The existence (and classification) on general
Poisson manifolds remained a hard problem for quite some time. Beside
the remarkable case of linear Poisson structures \cite{gutt:1983a} not
much progress was achieved until the seminal work of Kontsevich
proving his formality conjecture from \cite{kontsevich:1997a} in the
preprint \cite{kontsevich:1997:pre}, later published in
\cite{kontsevich:2003a}. Tamarkin gave another interpretation of
Kontsevich's formality theorem using the language of operads in
\cite{tamarkin:1998b}. A different approach to globalization of
Kontsevich's local explicit formula is due to Dolgushev
\cite{dolgushev:2005a, dolgushev:2005b}. An interpretation of the
formality in terms of the Poisson sigma-model
\cite{schaller.strobl:1994a} can be found in the work of Cattaneo,
Felder and Tomassini \cite{cattaneo.felder.tomassini:2002b,
  cattaneo.felder.tomassini:2002a, cattaneo.felder:2001a,
  cattaneo.felder:2001d, cattaneo.felder:2000a}.

In addition to the existence and classification results many more
details on deformation quantization have been discussed over the
years. To mention just a few achievements we would like to point out
the construction of star products on Kähler manifolds taking advantage
of the additional complex structure. Here Cahen, Gutt, and Rawnsley
showed in \cite{cahen.gutt.rawnsley:1995a, cahen.gutt.rawnsley:1994a,
  cahen.gutt.rawnsley:1993a, cahen.gutt.rawnsley:1990a} how the
quantization scheme based on Berezin-Toeplitz operators
\cite{berezin:1975a, berezin:1975b, berezin:1975c} can be used to
obtain formal star products by asymptotic expansions of the integrals
encoding the Berezin-Toeplitz operators. In
\cite{bordemann.meinrenken.schlichenmaier:1991a} this was refined to
incorporate detailed operator norm estimates ultimately resulting in a
continuous field of $C^*$-algebras. The star products obtained this
way turn out to be of Wick type in the sense that one function is
differentiated in holomorphic directions while the other is
differentiated in anti-holomorphic directions only. Such star products
have been shown to exist in general by Karabegov
\cite{karabegov:1996a}, see also \cite{bordemann.waldmann:1997a} for
an explicit construction based on Fedosov's approach. Beside Kähler
manifolds, cotangent bundles are perhaps the most important classical
phase spaces. Deformation quantizations adapted to this class of
examples were discussed in \cite{pflaum:1998c, pflaum:1998b,
  bordemann.neumaier.pflaum.waldmann:2003a,
  bordemann.neumaier.waldmann:1998a,
  bordemann.neumaier.waldmann:1999a}.

A large activity over the years was on the understanding of star
products with symmetries: here one is given a classical symmetry in
form of a Lie group action or a Lie algebra action compatible with the
classical Poisson structure. Then the question is whether there is a
star product such that the Lie group acts by automorphisms of it, or,
infinitesimally, the Lie algebra acts by derivations. It turns out
that the existence of an invariant covariant derivative is sufficient
for this in general, see e.g. \cite{arnal.cortet.molin.pinczon:1983a,
  bertelson.bieliavsky.gutt:1998a, dolgushev:2005a,
  fedosov:1996a}. Classically, a symmetry becomes most interesting if
it is implemented by means of a momentum map. Again, the question is
whether one also has a corresponding momentum map on the quantum
side. Here one has several positive answers
\cite{mueller-bahns.neumaier:2004a, mueller-bahns.neumaier:2004b,
  gutt.rawnsley:2003a} including a complete classification in the
symplectic case \cite{reichert.waldmann:2016a}. One important
construction for classical mechanical systems with symmetry is the
Marsden-Weinstein reduction \cite{marsden.weinstein:1974a} which
allows to reduce the dimension by fixing the values of conserved
quantities build out of the momentum map. In deformation quantization
several quantum analogs of reduction are constructed by Fedosov
\cite{fedosov:1998a}, by Bordemann, Herbig, and Waldmann
\cite{bordemann.herbig.waldmann:2000a} based on a BRST formulation as
well as by Gutt and Waldmann \cite{gutt.waldmann:2010a} including a
careful discussion of the involutions. Recently, Reichert computed the
characteristic classes of reduced star products in terms of the
equivariant classes on the original phase space in
\cite{reichert:2017a}. On the more conceptual side, Cattaneo and
Felder investigated the compatibility of formality morphisms with
coisotropic submanifolds in general \cite{cattaneo.felder:2007a,
  cattaneo.felder:2004a, cattaneo:2004a}.

Finally, a last important aspect of deformation quantization is the
development of a physically reasonable notion of states, considered as
positive linear functionals on the algebra of observables
\cite{bordemann.waldmann:1998a}. Here positivity is understood in the
sense of the canonical ring ordering of $\mathbb{R}[[\lambda]]
\subseteq \mathbb{C}[[\lambda]]$. In fact, large parts of operator
algebraic theory on states and representations can be transferred to
this purely algebraic framework based on $^*$-algebras over ordered
rings, see e.g. \cite{bordemann.roemer.waldmann:1998a,
  bursztyn.waldmann:2001b, bursztyn.waldmann:2001a, waldmann:2005b} as
well as the book \cite[Chap.~7]{waldmann:2007a} and the lecture notes
\cite{waldmann:2019a:script}. Beside studying the existence of
positive functionals as deformations of classical ones
\cite{bursztyn.waldmann:2005a, bursztyn.waldmann:2000a} one can
establish a notion of strong Morita equivalence yielding the
equivalence of categories of $^*$-representations
\cite{bursztyn.waldmann:2000b, bursztyn:2001b, bursztyn:2001a,
  bursztyn.waldmann:2001a, bursztyn:2002a, bursztyn.waldmann:2004a,
  bursztyn.waldmann:2005c, bursztyn.waldmann:2005b,
  bursztyn.waldmann:2012b, bursztyn.waldmann:2012a} culminating in a
complete and geometrically simple classification of star products up
to Morita equivalence in \cite{ bursztyn.waldmann:2002a,
  bursztyn.dolgushev.waldmann:2012a}, first in the symplectic and then
in the general Poisson case. Conversely, the investigations of the
Picard group of star products triggered an analogous program also for
the semi-classical counterpart, the Picard group(oid) in Poisson
geometry \cite{bursztyn.weinstein:2005a, bursztyn.weinstein:2004a,
  bursztyn.fernandes:2018a}. Again, one has adapted versions
including classical symmetries \cite{jansen.waldmann:2006a,
  jansen.neumaier.schaumann.waldmann:2012a}.

Finally, it should be mentioned that star products found their way
also to field-theoretic applications: while the original intention was
to provide a quantization scheme for geometrically non-trivial but
finite-dimensional phase spaces the notion of star products and hence
deformation quantization is applicable also to classical systems with
infinitely many degrees of freedom, i.e. classical field theories. The
main difficulty is to find a suitable Hamiltonian formulation first
with a Poisson algebra of observables. Once this is achieved, the very
definition of a star product clearly makes sense and can be
explored. In \cite{dito:2002a, dito:1993a, dito:1992a, dito:1990a}
Dito investigated this possibility for the Klein-Gordon field theory
and studied the effects of renormalization needed for interacting
fields. Fredenhagen and Dütsch continued and extended this
considerably \cite{duetsch.fredenhagen:1999a,
  duetsch.fredenhagen:2001b, duetsch.fredenhagen:2001a,
  duetsch.fredenhagen:2004a, duetsch.fredenhagen:2003a} and ever since
star products became an important tool in more conceptual approaches
to quantum field theory on globally hyperbolic spacetimes, see
e.g. \cite{fredenhagen.rejzner:2016a} as well as the monograph
\cite{rejzner:2016a} for further references.

%
% The Quest for Convergence
%

\section{The Quest for Convergence}
\label{sec:QuestConvergence}

All the above achievements will not allow to hide the open question on
the convergence of the formal star products. It is clear that for a
physically applicable quantization one can not treat $\hbar$ as a
formal parameter but needs to paste in a numerical value determined in
a given system of units. If not on the level of observable algebras,
the latest moment where this definitely has to be done is for
expectation values and spectral values: ultimately, they are
measurable quantities accessible by experiments and here we have no
formal series anymore.

The problem becomes now manifest since in all known examples of star
products the operators $C_r$ are bidifferential operators of order at
least $r$ in each argument. In particular, the star product
$f \star g$ of two smooth functions will typically see and use
essentially the whole infinite jet information of $f$ and $g$. This
shows that one can find immediately two functions
$f, g \in \Cinfty(M)$ in such a way that their (formal) star product
$f \star g$ at a given point $p \in M$ will have radius of convergence
in $\lambda$ being zero. The reasons is simply that with the classical
Borel lemma we can adjust the Taylor coefficients of $f$ and $g$ at a
point $p$ arbitrarily bad, leading to a divergent series $f \star g$
unless $\lambda = 0$. Thus the class of \emph{all} smooth functions is
not suitable for a naive (pointwise) convergence of $f \star g$.
There are now several ways considered to circumvent this fundamental
problem.

The perhaps first interpretation of the convergence problem is that
the formal star product $f \star g$ should be interpreted as an
\emph{asymptotic expansion} of some honest product $f \circ_\hbar g$
for $\hbar \longrightarrow 0^+$, at least for a suitably interesting
class of functions. In fact, this seems not to be completely hopeless
as in many examples the formal star products where exactly constructed
that way. One investigated certain integral formulas for operators and
compositions of integral operators depending on a positive parameter
$\hbar$ and showed that in the limit $\hbar \longrightarrow 0$ one
obtains a smooth dependence on $\hbar$ (from the right). Thus the
resulting formal Taylor expansion in $\hbar$ of these integral
expressions yielded a composition law $f \star g$ which then turned
out to be a formal star product.

This has been a successful approach in several situations, most
notably in the early constructions of star products as in
\cite{cahen.gutt.rawnsley:1995a, cahen.gutt.rawnsley:1994a,
  cahen.gutt.rawnsley:1993a, cahen.gutt.rawnsley:1990a,
  bordemann.meinrenken.schlichenmaier:1991a}. The tricky question in
these approaches is, however, whether one actually has an honest
\emph{algebra} of functions, i.e. whether one can specify a subspace
of functions which is indeed closed under the integral formulas for
the product. One might be able to circumvent this by starting with a
small class of functions for which the integral formulas are easily
defined, say with suitable support conditions, and define the algebra
implicitly as all what is generated by the products of such nice
functions. It remains a nontrivial question whether this is really a
well-defined algebra as the integral formulas might produce functions
for which a further application of the integration is no longer valid.

On a more conceptual level, the very few integral formulas which are
known are used to obtain universal deformation formulas. The most
prominent example is the integral formula for the Weyl-Moyal star
product. It can be viewed as an integral kernel on the abelian Lie
group $\mathbb{R}^{2n}$. Whenever one has now an action of this group
on a suitable topological algebra, say a Banach algebra or a locally
convex algebra, subject to certain continuity requirements as
e.g. strongly continuous and polynomially bounded, then one can use
the action to induce a deformation of the algebra. This way one
obtains a new product, still being continuous, which depends on
$\hbar$ in a very controlled way. The indicated procedure can be
traced back to Rieffel's seminal work \cite{rieffel:1993a}, see also
his works \cite{rieffel:1989b, rieffel:1989a} where the main focus is
on the case of $C^*$-algebras. Here one wants to construct in addition
a $C^*$-norm on the deformed algebra such that it becomes a
$C^*$-algebra after completion. Ultimately, this results in a
continuous field of $C^*$-algebras parametrized by $\hbar \in [0,
\infty)$. Many more details on this construction and its application
in quantization theory can be found in Landsman's monograph
\cite{landsman:1998a}, a generalization to general locally convex
algebras and their modules under very mild assumptions on the group
action is discussed in \cite{lechner.waldmann:2016a}.

However, the above class of constructions only works if the abelian
group $\mathbb{R}^{2n}$ acts. There are many interesting situations in
symplectic geometry where one immediately finds obstructions that the
Poisson structure comes from an action of $\mathbb{R}^{2n}$, or, even
worse, from any kind of Lie group action, see
\cite{bieliavsky.esposito.waldmann.weber:2018a}. It took some
considerable effort to find integral formulas beyond the abelian
case. Here Bieliavsky and Gayral found a vast generalization of
Rieffel's original ideas and proposed universal deformation formulas
in \cite{bieliavsky.gayral:2015a}, see also the earlier works
\cite{bieliavsky.detournay.spindel:2009a, bieliavsky:2008a,
  bieliavsky.bonneau.maeda:2007a, bieliavsky:2002a,
  bieliavsky.massar:2001b, bieliavsky.massar:2001a}. These ideas lead
ultimately to a quantization of Riemann surfaces of higher genus
\cite{bieliavsky:2017a}. Yet a different approach to quantum surfaces
in a $C^*$-algebraic formulation is due to Natsume, Nest, and Peters
\cite{natsume.nest.peter:2003a, natsume:2000a, natsume.nest:1999a}
providing strict quantizations, too.

It is now a second option we would like to advocate for: instead of
replacing the formal power series by integral formulas and
re-interpreting the star product in such terms, we would like to take
the series description serious and investigate the actual convergence
of the formal power series in $\hbar$ for certain classes of
functions. To achieve this goal, one has to master several steps:
\begin{enumerate}
\item First a class of function has to be determined on which the
    formal star product is known to converge. In many examples this is
    not too bad and one has good candidates. In fact, in the examples
    we discuss below this class of functions is always more or less
    related to polynomials and thus the formal power series simply
    will \emph{terminate} after finitely many contributions.  However,
    this naive class of functions is typically very small and too
    non-interesting for serious applications.
\item Thus the second step consists in finding a suitable locally
    convex topology on the above naive class of functions for which
    the product is continuous. This is typically a real extra piece of
    information and it seems that there is no canonical way to achieve
    this. In the examples one has several possibilities and it is not
    completely clear how to characterize an optimal solution among
    many.
\item The third step then consists in extending the product by
    continuity to the completion of the naive class of functions. Here
    one has hopefully found a coarse enough topology such that the
    completion is interesting enough.
\item In a last step one can now try to examine the completion and
    determine whether its elements can still be considered as
    functions on the original phase space. This essentially means that
    the evaluation functionals at points $p \in M$ are continuous with
    respect to the found locally convex topology. In addition one
    wants to determine things like positive functionals etc.
\end{enumerate}

While there seems to be no general theory in sight which would allow
to make conceptual statements about this program, it nevertheless
turns out to work well in several important classes of examples. What
is perhaps more important is the perspective towards
infinite-dimensional situations. Here the integral formulas clearly
stop to make sense and one needs a replacement. The question about
convergence, however, still can be asked and, in many examples, be
answered to the positive. Thus the above program might not only be
interesting within the realm of deformation quantization of classical
mechanical systems but also beyond in field-theoretic situations.

%
% The Weyl Star Product
%

\section{The Weyl Star Product}
\label{sec:WeylStarProduct}

We start with the perhaps most important example in deformation
quantization, the Weyl star product and its relatives. These are star
products quantizing the constant (symplectic) Poisson structure on a
vector space. In this section we follow closely the construction from
\cite{waldmann:2014a} and \cite{schoetz.waldmann:2018a}. The previous
constructions from \cite{beiser.roemer.waldmann:2007a} only work in
the finite-dimensional case and for a particular case, the Wick star
product, where they produce a slightly finer topology.

To set the stage we consider a real vector space $V$ with a bilinear
form $\Lambda\colon V \times V \longrightarrow \mathbb{C}$. In many
cases the values are real but we leave the option to have complex
values at the moment. The Weyl product will then be defined on
polynomials. Now in finite dimensions, the real-valued polynomials on
the dual space $V^*$ are just the symmetric algebra
$\Sym_{\mathbb{R}}^\bullet(V)$ while for infinite-dimensional vector
spaces this is of course no longer true. For technical reasons it is
convenient to work with the symmetric algebra
$\Sym_{\mathbb{R}}^\bullet(V)$ over $V$ as a replacement. They can be
viewed as the polynomials on the pre-dual of $V$ if such a pre-dual
exists: we do not require this but consider
$\Sym_{\mathbb{R}}^\bullet(V)$ only.

We then consider the complexified symmetric algebra
$\Sym^\bullet(V) = \Sym^\bullet_{\mathbb{C}}(V)$ of $V$ as replacement
for the complex-valued polynomials. Out of $\Lambda$ we obtain the
operator
\begin{equation}
    \label{eq:PoissonOperatorPLambda}
    P_\Lambda\colon
    \Sym^\bullet(V) \tensor \Sym^\bullet(V)
    \longrightarrow
    \Sym^\bullet(V) \tensor \Sym^\bullet(V)
\end{equation}
defined on factorizing symmetric tensors
$v_1 \cdots v_n \in \Sym^n(V)$ and $w_1 \cdots w_m \in \Sym^m(V)$ by
\begin{equation}
    \label{eq:PLambdaDef}
    P_\Lambda(v_1 \cdots v_n \tensor w_1 \cdots w_m)
    =
    \sum_{k = 1}^n \sum_{\ell = 1}^m
    \Lambda(v_k, w_\ell)
    v_1 \cdots \stackrel{k}{\wedge} \cdots v_n
    \tensor
    w_1 \cdots \stackrel{\ell}{\wedge} \cdots w_m
\end{equation}
and extended linearly. Here $\stackrel{k}{\wedge}$ means to omit the
$k$-th factor and the empty products are defined to be
$\Unit \in \Sym^0(V) = \mathbb{C}$ and finally we set
$P_\Lambda(\Unit \tensor w) = 0 = P_\Lambda(v \tensor \Unit)$ for all
$v, w \in \Sym^\bullet(V)$.
\begin{lemma}
    \label{lemma:ConstantPoisson}%
    The operator $P_\Lambda$ satisfies a Leibniz rule in each tensor
    factor and gives a Poisson structure by
    \begin{equation}
        \label{eq:ConstantPoisson}
        \{v, w\}
        =
        \mu \circ \big(
        P_\Lambda - \tau \circ P_\Lambda \circ \tau
        \big)
        (v \tensor w)
    \end{equation}
    for $v, w \in \Sym^\bullet(V)$.
\end{lemma}
Here $\mu$ denotes the symmetric tensor product and $\tau$ is the
canonical flip of two tensor factors. The proof consists in
considering the following operators
\begin{equation}
    \label{eq:Peinszweidrei}
    P_{12}, P_{13}, P_{23}\colon
    \Sym^\bullet(V) \tensor \Sym^\bullet(V) \tensor \Sym^\bullet(V)
    \longrightarrow
    \Sym^\bullet(V) \tensor \Sym^\bullet(V) \tensor \Sym^\bullet(V)
\end{equation}
defined by $P_{12} = P_\Lambda \tensor \id$,
$P_{23} = \id \tensor P_\Lambda$ and
$P_{13} = (\id \tensor \tau) \circ (P_\Lambda \tensor \id) \circ (\id
\tensor \tau)$.
It can now be shown by a straightforward computation that they
pairwise commute and satisfy the Leibniz rules
\begin{equation}
    \label{eq:PeinszweidreiRelations}
    P_\Lambda \circ (\mu \tensor \id)
    =
    (\mu \tensor \id) \circ (P_{13} + P_{23})
    \quad
    \textrm{and}
    \quad
    P_\Lambda \circ (\id \tensor \mu)
    =
    (\id \tensor \mu) \circ (P_{13} + P_{12}).
\end{equation}
From this the properties of a Poisson bracket follow
immediately. Moreover, they allow to obtain an associative deformation
at once. Note that the operator $P_\Lambda$ lowers the symmetric
degrees in both tensor factors by one and hence the following
exponential series
\begin{equation}
    \label{eq:WeylStarProduct}
    v \star_{z\Lambda} w
    =
    \mu \circ \E^{z P_\Lambda} (v \tensor w)
\end{equation}
is well-defined for $v, w \in \Sym^\bullet(V)$ and for all
$z \in \mathbb{C}$. This is the Weyl star product:
\begin{lemma}
    \label{lemma:WeylStarProduct}%
    For all $z \in \mathbb{C}$ the product $\star_{z\Lambda}$ is an
    associative product for $\Sym^\bullet(V)$. In zeroth order of $z$,
    i.e. for $z = 0$, it yields the symmetric tensor product and the
    commutator in first order of $z$ is the Poisson bracket
    \eqref{eq:ConstantPoisson}.
\end{lemma}
Strictly speaking, the Weyl product would be obtained from an
antisymmetric $\Lambda$ taking only real values and $z = \I\hbar$.
Allowing for symmetric contributions of $\Lambda$ and possible real
values will give other star products quantizing the constant Poisson
structure like a Wick star product or a standard-ordered star product.
Nevertheless, we will simply speak of the Weyl product in the
following.

The important point is that the product $\star_{z\Lambda}$ converges
for trivial reasons on the polynomials, i.e. on $\Sym^\bullet(V)$,
since the exponential series simply terminates after finitely many
contributions.

It is now the second step of the program which requires some more
effort. We need to find a suitable topology on $\Sym^\bullet(V)$ which
will turn $\star_{z\Lambda}$ into a continuous product.

Here we have to require a locally convex topology for $V$ and
continuity properties for $\Lambda$. In the finite-dimensional case
this is automatic as then $V$ has a unique Hausdorff locally convex
topology (even coming from Hilbert space structures) and every
bilinear form is continuous. In the infinite-dimensional case this is
an additional information we have to invest.  Thus assume $V$ is
locally convex and $\Lambda$ is continuous: we require continuity and
not just separate continuity.

For each tensor power $\Tensor^k(V)$ and hence for each symmetric
power $\Sym^k(V) \subseteq \Tensor^k(V)$ of $V$ we can then use the
projective topology: this is the locally convex topology determined by
the seminorms
\begin{equation}
    \label{eq:pnDef}
    \halbnorm{p}^n = \halbnorm{p} \tensor \cdots \tensor \halbnorm{p}
\end{equation}
for all continuous seminorms $\halbnorm{p}$ of $V$. As usual, it
suffices to consider a defining system of continuous seminorms
$\halbnorm{p}$ of $V$ to obtain a defining system $\halbnorm{p}^n$ for
the projective topology of $\Tensor^k(V)$ and hence for $\Sym^k(V)$.
For $n = 0$ we define $\halbnorm{p}^0$ to be the usual absolute value
on $\mathbb{R}$ and $\mathbb{C}$.

To generate a topology on the tensor algebra and the symmetric algebra
is now less canonical. Here we fix a parameter $R \in \mathbb{R}$ to
parametrize different possibilities for locally convex topologies
\cite[Def.~3.5]{waldmann:2014a}:
\begin{definition}[$\Tensor_R$- and $\Sym_R$-topology]
    \label{definition:TRtopology}%
    Let $R \in \mathbb{R}$ and let $V$ be a locally convex vector
    space. On the tensor algebra $\Tensor^\bullet(V)$ of $V$ one
    defines the seminorm
    \begin{equation}
        \label{eq:pRseminorm}
        \halbnorm{p}_R(v)
        =
        \sum_{n=0}^\infty
        n!^R \halbnorm{p}^n(v_n)
    \end{equation}
    for a continuous seminorm $\halbnorm{p}$ on $V$ where
    $v = \sum_{n \in \mathbb{N}_0} v_n$ are the homogeneous components
    of $v \in \Tensor^\bullet(V)$. The locally convex topology
    generated by all these seminorms on $\Tensor^\bullet(V)$ for
    $\halbnorm{p}$ varying through all continuous seminorms on $V$ is
    called the $\Tensor_R$-topology. The subspace topology induced on
    $\Sym^\bullet(V)$ is called the $\Sym_R$-topology.
\end{definition}
The tensor algebra (symmetric algebra) equipped with the
$\Tensor_R$-topology ($\Sym_R$-topology) will be denoted by
$\Tensor^\bullet_R(V)$ and $\Sym^\bullet_R(V)$, respectively.
\begin{remark}[Properties of the $\Tensor_R$- and $\Sym_R$-topology]
    \label{remark:PropertiesTR}%
    Let $V$ be a locally convex vector space. The following properties
    of the $\Tensor_R$- and $\Sym_R$-topology were obtained in
    \cite{waldmann:2014a}:
    \begin{remarklist}
    \item \label{item:TRInducesPi} For every $R \in \mathbb{R}$ and
        every $k \ge 0$ the induced topology on
        $\Tensor^k(V) \subseteq \Tensor^\bullet(V)$ is the projective
        topology. The same holds for the symmetric
        version. Conversely, the inclusion
        $\Tensor^k(V) \subseteq \Tensor^\bullet(V)$ is continuous with
        respect to the projective topology and the
        $\Tensor_R$-topology, respectively.
    \item \label{item:TRHausdorff} For every $R \in \mathbb{R}$ the
        $\Tensor_R$-topology is Hausdorff iff $V$ is Hausdorff. In
        this case also the $\Sym_R$-topology is Hausdorff.
    \item \label{item:EquivalentViaSup} An equivalent system of
        seminorms for the $\Tensor_R$-topology is obtained by taking
        all the seminorms
        \begin{equation}
            \label{eq:pRinfty}
            \halbnorm{p}_{R, \infty}(v)
            =
            \sup_{n \in \mathbb{N}_0} n!^R \halbnorm{p}^n(v_n)
        \end{equation}
        for all continuous seminorms $\halbnorm{p}$ on $V$.
    \item \label{item:RRprimeFiner} Let $R' > R$. Then the
        $\Tensor_R$-topology is coarser than the
        $\Tensor_{R'}$-topology.
    \item \label{item:CompletionTR} The completion of
        $\Tensor^\bullet(V)$ with respect to the $\Tensor_R$ topology
        is explicitly given by
        \begin{equation}
            \label{eq:CompletionTR}
            \complete{\Tensor}_R(V)
            =
            \bigg\{
                v = \sum_{n=0}^\infty v_n
                \; \bigg| \;
                \sum_{n=0}^\infty n!^R \halbnorm{p}^n(v_n) < \infty
                \textrm{ for all }
                \halbnorm{p}
            \bigg\}
            \subseteq
            \prod_{n=0}^\infty V^{\hat{\tensor} n},
        \end{equation}
        where $V^{\hat{\tensor} n}$ denotes the completion of the
        $n$-th tensor power in the projective topology. Analogously,
        one obtains the completion $\complete{\Sym}^\bullet_R(V)$ of
        $\Sym^\bullet_R(V)$ as those series in
        $\complete{\Tensor}^\bullet_R(V)$ consisting of symmetric
        tensors in each degree $n \in \mathbb{N}_0$.
    \item \label{item:LCAlgebras} The tensor product as well as the
        symmetric tensor product are continuous products for
        $R \ge 0$. For $R = 0$ the tensor algebra
        $\Tensor^\bullet_0(V)$ is the free locally multiplicatively
        convex algebra generated by $V$ and $\Sym_0(V)$ is the free
        locally multiplicatively convex commutative algebra generated
        by $V$ as used e.g. by Cuntz in \cite{cuntz:1997a}. For
        $R > 0$ the topologies are \emph{not} locally multiplicatively
        for the (symmetric) tensor product anymore.
    \item \label{item:FirstCountable} The $\Tensor_R$- and the
        $\Sym_R$-topology are first countable iff the original
        topology of $V$ is first countable.
    \item \label{item:ContinuousFunctionals} Let $R \ge 0$. Then for
        every continuous linear functional $\varphi \in V'$ on $V$,
        the corresponding evaluation functional
        \begin{equation}
            \label{eq:CharacterContinuous}
            \delta_\varphi\colon \Tensor^\bullet_R(V)
            \ni
            v = \sum_{n=0}^\infty v_n
            \; \mapsto \;
            \delta_\varphi(v)
            =
            \sum_{n=0}^\infty \varphi^{\tensor n}(v_n)
            \in \mathbb{C}
        \end{equation}
        is a continuous algebra homomorphism with respect to the
        tensor product. In particular,
        $\delta_\varphi\colon \Sym^\bullet_R(V) \longrightarrow
        \mathbb{C}$
        is a character and hence we can identify $\Sym^\bullet(V)$
        with certain polynomial functions on $V'$.
    \end{remarklist}
\end{remark}
All these properties are rather straightforward to check. Note that
the statement about the completion shows that we indeed get an
interesting completion the smaller $R$ becomes. Already for $R < 1$ we
have e.g. the exponential series $\exp(v)$ for $v \in V$ in the
completion.

Slightly more important are the next two properties which show that
these topologies preserve interesting properties of the underlying
vector space $V$. Note that in the finite-dimensional case we are in
the situation to apply the following results
\cite[Thm.~4.10]{waldmann:2014a}:
\begin{theorem}
    \label{theorem:NuclearSchauder}%
    Let $V$ be a Hausdorff locally convex space and let $R \ge 0$.
    \begin{theoremlist}
    \item \label{item:Schauder} The space $V$ admits an absolute
        Schauder basis iff $\Tensor_R(V)$ admits an absolute Schauder
        basis iff $\Sym_R(V)$ admits an absolute Schauder basis.
    \item \label{item:Nuclear} The space $V$ is nuclear iff
        $\Tensor_R(V)$ is nuclear iff $\Sym_R(V)$ is nuclear.
    \end{theoremlist}
\end{theorem}

We come now to the main result of this section: the continuity of the
Weyl product. Here we have the following statement first formulated in
this generality in \cite[Thm.~3.17]{waldmann:2014a}:
\begin{theorem}[Weyl product]
    \label{theorem:ContinuityWeyl}%
    Let $R \ge \frac{1}{2}$. Then the Weyl product $\star_{z\Lambda}$
    is continuous with respect to the $\Sym_R$-topology. Moreover, for
    all elements $v, w \in \complete{\Sym}_R(V)$ in the completion the
    series
    \begin{equation}
        v \star_{z\Lambda} w
        =
        \sum_{n=0}^\infty \frac{z^n}{n!}
        \mu(P_{\Lambda}^n (v \tensor w))
    \end{equation}
    converges absolutely in the $\Sym_R$-topology. The dependence on
    $z \in \mathbb{C}$ is entire.
\end{theorem}
This settles the questions raised in
Section~\ref{sec:QuestConvergence} for the Weyl product
$\star_{z\Lambda}$ completely. Note that the continuity of the
characters allows to interpret the elements in the completion
$\complete{\Sym}_R^\bullet(V)$. We conclude this section now with a
few remarks:
\begin{remark}[Grassmann version]
    \label{remark:GrassmannVersion}%
    In \cite{waldmann:2014a} a slightly more general situation was
    considered: the vector space $V$ was allowed to carry a
    $\mathbb{Z}_2$-grading which then was used to incorporate
    additional signs in the definition of the symmetric tensor product
    and $P_\Lambda$. Then the graded symmetric algebra is the usual
    symmetric algebra for even vectors but the Grassmann algebra for
    odd vectors. In total one can obtain a combination of both,
    thereby allowing also a deformation quantization of both. The even
    part is then quantized by the above Weyl product, for the
    Grassmann part one obtains a quantization by a Clifford
    algebra. Needless to say that the analysis and continuity
    estimates will ignore the signs and are thus valid also in this
    slightly more general situation.
\end{remark}
\begin{remark}[Finite dimensions]
    \label{remark:FiniteDim}%
    In finite dimensions \cite{omori.maeda.miyazaki.yoshioka:2002a,
      omori.maeda.miyazaki.yoshioka:2000a,
      omori.maeda.miyazaki.yoshioka:2007a} considered already a
    particular case of the Weyl product and discussed the convergence
    issues. They found a classical functional space of entire
    functions subject to certain growth conditions at infinity for
    which the Weyl product converges. It turns out that this
    functional space coincides with the above completion
    $\complete{\Sym}_R(V)$ for an appropriate choice of $R$.
\end{remark}
\begin{remark}[Pro-Hilbert case]
    \label{remark:ProHilbertCase}%
    A more recent development is proposed in
    \cite{schoetz.waldmann:2018a} where the locally convex space $V$
    is considered to be a projective limit of Hilbert spaces. This
    includes the case where $V$ is actually a Hilbert space. But also
    every nuclear space is projectively Hilbert. In this situation one
    can use the inner products defining the Hilbert seminorms to
    establish a different topology on each tensor power
    $V^{\tensor n}$: one simply extends the inner products and
    computes the corresponding seminorm afterwards. This results in a
    slightly coarser topology than the projective topology on
    $V^{\tensor n}$. In some sense this corresponds more to a
    Hilbert-Schmidt like topology rather than a trace-class
    topology. Now the completion becomes more interesting already for
    fixed $n$. The continuity of the Weyl product can now be shown
    also in this case with an extension of the topology to the whole
    symmetric algebra as before: in \cite{schoetz.waldmann:2018a} we
    focused on the case $R = \frac{1}{2}$ directly, yielding the
    coarsest possible topology. Many additional feature including a
    detailed description of the Gel'fand transform in the classical
    case can be done. This allows to determine the completion very
    explicitly as certain real-analytic functions on the topological
    dual $V'$. Finally, in the nuclear case we can benefit from both
    versions as here the two competing topologies on the tensor powers
    simply coincide. Thus all properties obtained above as well as
    those from \cite{schoetz.waldmann:2018a} become
    available. Luckily, the physically relevant cases seem all to be
    of that type: finite-dimensional spaces as well as the usual test
    function spaces used in (quantum) field theory are all nuclear.
\end{remark}

%
% The Gutt Star Product
%

\section{The Gutt Star Product}
\label{sec:GuttStarProduct}

As a next important class of examples one can consider the linear
Poisson structures on the dual of a Lie algebra. This has been
investigated in Poisson geometry from many points of view and provides
one of the most important examples as it shows that Lie algebra theory
becomes accessible via Poisson geometric techniques. In this section
we discuss the continuity and convergence results form
\cite{esposito.stapor.waldmann:2017a}.

To define the linear Poisson structure we consider a real Lie algebra
$\lie{g}$ and instead of complex-valued polynomial functions on its
dual $\lie{g}^*$ we focus on the symmetric algebra
$\Sym^\bullet(\lie{g})$ for the same reasons as already in the Weyl
product case. In particular, we allow for infinite-dimensional Lie
algebras as well.

The star product we want to consider originates from the
Poincaré-Birkhoff-Witt isomorphism between the (complexified)
symmetric algebra and the (complexified) universal enveloping algebra
$\mathcal{U}(\lie{g})$. In some more detail one considers the $n$-th
symmetrization map defined on factorizing symmetric tensors by
\begin{equation}
    \label{eq:PBWIso}
    \mathfrak{q}_n\colon \Sym^n(\lie{g})
    \ni
    \xi_1 \cdots \xi_n
    \; \mapsto \;
    \frac{1}{n!} \sum_{\sigma \in S_n}
    \xi_{\sigma(1)} \odot \cdots \odot \xi_{\sigma(n)}
    \in
    \mathcal{U}(\lie{g}),
\end{equation}
and extended linearly. Here $\odot$ denotes the product of the
universal enveloping algebra. Then the Poincaré-Birkhoff-Witt theorem
says that the direct sum
$\mathfrak{q} = \sum_{n=0}^\infty \mathfrak{q}_n$ provides a vector
space isomorphism between $\Sym^\bullet(\lie{g})$ and
$\mathcal{U}(\lie{g})$.

Denote the projection onto the $n$-th homogeneous component in
$\Sym^\bullet(\lie{g})$ by
$\pr_n\colon \Sym^\bullet(\lie{g}) \longrightarrow \Sym^n(\lie{g})$.
Then for $z \in \mathbb{C}$ and homogeneous elements
$x \in \Sym^k(\lie{g})$ and $y \in \Sym^\ell(\lie{g})$ one defines the
product
\begin{equation}
    \label{eq:GuttHomogeneous}
    x \star_z y
    =
    \sum_{n = 0}^{k + \ell - 1}
    z^n \pr_{k+\ell-n}
    \big(
    \mathfrak{q}^{-1}(\mathfrak{q}(x) \odot \mathfrak{q}(y))
    \big)
\end{equation}
and extends this bilinearly to $\Sym^\bullet(\lie{g})$. Again, we see
that this is a well-defined product for all $z \in \mathbb{C}$ which
turns out to be associative. In fact, the case $z = 1$ is clear as
there it becomes isomorphic to the universal enveloping algebra
directly. The other values of $z$ can be understood as a rescaling of
the Lie bracket by $z$. In particular, for $z = 0$ one is back at the
symmetric algebra. The first order term in $z$ is then the usual
linear Poisson bracket for $\Sym^\bullet(\lie{g})$ obtained from the
Lie bracket of $\lie{g}$. This can equivalently be obtained by
extending $[\argument, \argument]_{\lie{g}}$ to
$\Sym^\bullet(\lie{g})$ enforcing the Leibniz rule as usual.
\begin{definition}[Gutt star product]
    \label{definition:GuttStarProduct}%
    The star product $\star_z$ is called the Gutt star product for
    $\lie{g}$.
\end{definition}
In fact, Gutt introduced this star product in \cite{gutt:1983a} as an
intermediate step to obtain a star product for the cotangent bundle of
a Lie group $G$ integrating $\lie{g}$. Alternatively and
independently, Drinfel'd introduced this star product in
\cite{drinfeld:1983a} in the context of quantum group theory.

The first important observation is the relation of the Gutt star
product with the Baker-Campbell-Hausdorff series. In fact, one has the
following statement:
\begin{lemma}
    \label{lemma:GuttBCH}%
    Viewing $\exp(\xi)$ as a formal series in
    $\prod_{n=0}^\infty \Sym^n(\lie{g})$ for $\xi \in \lie{g}$ the
    Gutt star product is determined by
    \begin{equation}
        \label{eq:expGuttBCH}
        \exp(\xi) \star_z \exp(\eta)
        =
        \exp\big(\tfrac{1}{z} \BCHseries(z\xi, z\eta)\big),
    \end{equation}
    where $\BCHseries$ is the usual Baker-Campbell-Hausdorff series of
    $\lie{g}$.
\end{lemma}
Here the idea is that by differentiation of the left hand side one can
obtain all monomials $\xi^k \star_z \eta^\ell$ with
$k, \ell \in \mathbb{N}_0$. Then the corresponding right hand side
gives an explicit formula for the Gutt star products of such
monomials. By polarization, this determines $\star_z$ also for
$\xi_1 \cdots \xi_k \in \Sym^k(\lie{g})$ and
$\eta_1 \cdots \eta_\ell \in \Sym^\ell(\lie{g})$ for arbitrary
$\xi_1, \ldots, \xi_k, \eta_1, \ldots, \eta_\ell \in \lie{g}$. Hence
the Gutt star product is completely encoded in the
Baker-Campbell-Hausdorff series of $\lie{g}$.
\begin{remark}
    \label{remark:IntegralFormulas}%
    The Gutt star product can also be formulated using integral
    formulas as this has (implicitly) been done already by Berezin
    \cite{berezin:1967a} and later on by Rieffel in
    \cite{rieffel:1990b}. However, for the analysis we have in mind
    the above version based on the Baker-Campbell-Hausdorff series is
    the most suitable. In particular, the integral formulas only make
    sense in finite dimensions while \eqref{eq:expGuttBCH} holds
    algebraically in general.
\end{remark}

Now we consider a locally convex topology on $\lie{g}$ in addition.
As in the constant case of the Weyl product we want to use the
$\Sym_R$-topology also for the linear case. And, again as in the
constant case, we need some continuity for the Lie bracket. It turns
out that mere continuity of the bilinear map
$[\argument, \argument]_{\lie{g}}\colon \lie{g} \times \lie{g}
\longrightarrow \lie{g}$
is unfortunately not sufficient. Instead, we have to require a
slightly stronger condition. To this end we recall the definition of
an \emph{asymptotic estimate} algebra, see
\cite[Def.~1.1]{esposito.stapor.waldmann:2017a} as well as
\cite{boseck.czichowski.rudolph:1981a}:
\begin{definition}[Asymptotic estimate algebra]
    \label{definition:AsymptoticEstimate}%
    Let $\algebra{A}$ be a Hausdorff locally convex algebra (not
    necessarily associative) with multiplication $\cdot$.
    \begin{definitionlist}
    \item \label{item:AE} A continuous seminorm $\halbnorm{q}$ is
        called an asymptotic estimate for the continuous seminorm
        $\halbnorm{p}$ if for all $n \in \mathbb{N}$ and all words
        $w_n(x_1, \ldots, x_n)$ of $n-1$ products of the $n$ elements
        $x_1, \ldots, x_n \in \algebra{A}$ with arbitrary position of
        placing brackets one has
        \begin{equation}
            \label{eq:AE}
            \halbnorm{p}(w_n(x_1, \ldots, x_n))
            \le
            \halbnorm{q}(x_1) \cdots \halbnorm{q}(x_n).
        \end{equation}
    \item \label{item:AEAlgebra} The algebra $\algebra{A}$ is called
        an asymptotic estimate algebra (AE-algebra) if every
        continuous seminorm has an asymptotic estimate.
    \end{definitionlist}
\end{definition}
Clearly, the case $n = 2$ already shows that the bilinear product is
continuous. Thus we have a stronger form of continuity in an
AE-algebra.
\begin{remark}[Finite-dimensional algebras are AE]
    \label{remark:LMC}%
    If $\algebra{A}$ is an associative AE-algebra then $\algebra{A}$
    becomes an AE-Lie algebra with respect to the commutator. Indeed,
    this is a trivial estimate. Moreover, if $\algebra{A}$ is locally
    multiplicatively convex then it is AE for trivial reasons. In
    particular, every finite-dimensional algebra is AE since it is
    locally multiplicatively convex. Thus the requirement is always
    fulfilled in finite dimensions but it becomes interesting in the
    infinite-dimensional case.
\end{remark}

Having now an AE-Lie algebra one arrives at the following continuity
statement for the Gutt star product, where we have to use the
$\Sym_R$-topology for $R \ge 1$ to obtain continuity
\cite[Thm.~1.2]{esposito.stapor.waldmann:2017a}:
\begin{theorem}[Gutt star product]
    \label{theorem:GuttStar}%
    Let $R \ge 1$ and let $\lie{g}$ be an AE-Lie algebra.
    \begin{theoremlist}
    \item \label{item:GuttContinuous} The Gutt star product $\star_z$
        is continuous with respect to the $\Sym_R$-topology for all
        $z \in \mathbb{C}$.
    \item \label{item:CompletionTopHopf} The completion
        $\complete{\Sym}_R(\lie{g})$ becomes a locally convex Hopf
        algebra with respect to $\star_z$ and the undeformed
        coproduct, antipode, and counit.
    \item \label{item:GuttHolomorphic} The map
        $\mathbb{C} \ni z \mapsto x \star_z y \in
        \complete{\Sym}_R(\lie{g})$
        is holomorphic for all $x, y \in \complete{\Sym}_R(\lie{g})$.
    \item \label{item:GuttFunctorial} The construction is functorial
        for continuous Lie algebra morphisms.
    \end{theoremlist}
\end{theorem}
Here we use the fact that the universal enveloping algebra
$\mathcal{U}(\lie{g})$ is a Hopf algebra in the usual
way. Transferring this Hopf algebra structure back to $\Sym(\lie{g})$
gives the usual cocommutative coproduct $\Delta$ with corresponding
antipode and counit being the projection $\epsilon = \pr_0$. In
particular, $\Delta$ does not depend on the Lie algebra structure at
all. For a locally convex Hopf algebra, the coproduct is then allowed
to take values in the completed tensor product instead of just the
algebraic tensor product.

The proof requires a rather technical and detailed analysis of the
Baker-Campbell-Hausdorff series. The idea is to find estimates for the
homogeneous parts $\BCHseries_n(\xi, \eta)$ of $\BCHseries(\xi, \eta)$ where we
have exactly $n$ letters $\xi$ or $\eta$ and hence $n-1$ Lie
brackets. The $n-1$ brackets can be estimated by the asymptotic
estimates \eqref{eq:AE} and one is left with the numerical constants
in front. Here a result of Goldberg \cite{goldberg:1956a} shows that
\begin{equation}
    \label{eq:BCHnEstimate}
    \BCHseries_n(\xi, \eta)
    =
    \sum_{
      \substack{
        \textrm{Lie words $w$ in} \\
        \textrm{$\xi$, $\eta$ of length $n$}
      }
    }
    \frac{g_w}{n} w
\end{equation}
with universal coefficients $g_w \in \mathbb{Q}$ satisfying the
estimate
\begin{equation}
    \label{eq:GoldbergEstimate}
    \sum_{\substack{\textrm{Lie words $w$} \\ \textrm{of length $n$}}}
    \frac{\abs{g_w}}{n}
    \le
    \frac{2}{n}.
\end{equation}
Using now the asymptotic estimate one finds for all Lie words $w$
containing $a$ times the letter $\xi$ and $b$ times the letter $\eta$
the estimate
\begin{equation}
    \label{eq:LieAEEstimate}
    \halbnorm{p}(w) \le \halbnorm{q}(\xi)^a\halbnorm{q}(\eta)^b.
\end{equation}
Form these two estimates one concludes now the continuity estimate for
the Gutt star product, once the correct combinatorics from
\eqref{eq:expGuttBCH} is taken into account, detail can be found in
\cite{esposito.stapor.waldmann:2017a}.

The result is somewhat disappointing as the condition $R \ge 1$ is
sharp: e.g. for the three-dimensional Heisenberg Lie algebra with
basis $Q, P, E$ and only nontrivial Lie bracket $[P, Q] = E$ one can
show that for $R < 1$ the product is discontinuous. The bound
$R \ge 1$ then forbids exponential series $\exp(\xi)$ to be in the
completion. Of course, this would have been too nice to have as then
the Lie group corresponding to $\lie{g}$ would have become part of the
completed algebra. On the other hand, this is of course not to be
expected as it would ultimately imply that the
Baker-Campbell-Hausdorff series would have radius of convergence being
infinite. And this one knows to be false in general.

Nevertheless, in the case of a \emph{nilpotent} Lie algebra one can
refine the above analysis in the following way. Instead of the
$\Sym_R$-topology for $R \ge 1$ one can consider the projective limit
of all $\Sym_R$-topologies for $R < 1$. It turns out that for a
locally convex nilpotent Lie algebra the Gutt star product is also
continuous with respect to this projective limit topology $\Sym_{1^-}$
and all the above results stay valid, see
\cite{esposito.stapor.waldmann:2017a}. Note that a nilpotent Lie
algebra is trivially AE as soon as the Lie bracket is continuous at
all.

%
% The Wick Type Star Product on the Poincaré Disc
%

\section{The Wick Type Star Product on the Poincaré Disc}
\label{sec:WickTypePoincareDisc}

The last example we want to discuss involves still a topologically
trivial phase space. However, now the underlying symplectic structure
comes from a curved Kähler structure. We consider the Poincaré disc
and its higher dimensional analogs. The star product of Wick type on
the Poincaré disc has a long history and was re-discovered many times,
see in particular the early contributions in
\cite{moreno.ortega-navarro:1983b, moreno:1986a,
  moreno.ortega-navarro:1983a} and \cite{cahen.gutt:1981a} as well as
the second part \cite{cahen.gutt.rawnsley:1993a}. Later on, the first
explicit formula was found in
\cite{bordemann.brischle.emmrich.waldmann:1996a,
  bordemann.brischle.emmrich.waldmann:1996b} together with some first
considerations on the convergence.  We will build on this construction
of the star product by a phase space reduction form $\mathbb{C}^{n+1}$
since this makes the already obtained convergence results for the Wick
star product available. This point of view was first taken in
\cite{beiser.waldmann:2014a} and recently substantially extended in
\cite{kraus.roth.schoetz.waldmann:2018a:pre}.

We first describe the construction of the Poincaré disc needed to
formulate the construction of the Wick star product for it. We
consider $\mathbb{C}^{n+1}$ with standard holomorphic coordinates
$(z^0, \ldots, z^n)$. We use an indefinite metric specified by the
matrix $g = \diag(-1, 1, \ldots, 1) \in \Mat_{n+1}(\mathbb{C})$ which
we identify with the corresponding quadratic function
\begin{equation}
    \label{eq:Functiong}
    g = g_{\mu\nu} z^\mu \cc{z}^\nu \in \Cinfty(\mathbb{C}^{n+1}).
\end{equation}
Here and in the following Greek indices run from $0$ to $n$ while
Latin indices will vary over $1$ to $n$ only. Moreover, we make use of
the summation convention throughout.  The function $g$ is then
invariant under the canonical and linear $\group{U}(1, n)$-action on
$\mathbb{C}^{n+1}$. This allows to consider the submanifold
\begin{equation}
    \label{eq:SubmanifoldZ}
    Z = g^{-1}(\{1\})
\end{equation}
of $\mathbb{C}^{n+1}$. Note that $1$ is indeed a regular value of $g$
and hence $Z$ is a submanifold. Moreover, $Z$ is clearly
$\group{U}(1, n)$-invariant.

In a next step we consider the central diagonal
$\group{U}(1) \subseteq \group{U}(1, n)$. Since this subgroup acts
freely (and properly for trivial reasons) on $Z$, the quotient
$D_n = Z/\group{U}(1)$ is a smooth manifold again and
$\pr\colon Z \longrightarrow D_n$ is a $\group{U}(1)$-principal fiber
bundle. The complementary subgroup
$\group{SU}(1, n) \subseteq \group{U}(1, n)$ still acts on $D_n$ and
provides a transitive smooth action. The manifold $D_n$ can now be
embedded into $\mathbb{CP}^n$ as follows. Note that $\group{U}(1, n)$
smoothly acts on $\mathbb{CP}^n$ by holomorphic diffeomorphisms with
three orbits, determined by the quadratic function $g$. The orbits are
characterized by the values of $g$: either positive or negative or $0$
on the complex lines in $\mathbb{C}^{n+1}$ representing the point in
$\mathbb{CP}^n$. Since $g$ is homogeneous, this is clearly
well-defined.  Then an equivalence class $\pr(p) \in D_n$ is mapped to
the equivalence class of the complex line through a representative
$p \in Z$. This gives a well-defined and $\group{U}(1, n)$-equivariant
embedding and identifies $D_n$ with the orbit of $\group{U}(1, n)$
where $g$ is positive. It follows that $D_n$ is an open subset of
$\mathbb{CP}^n$ and as such it inherits the structure of a complex
manifold.

In a last step we need the symplectic structure for $D_n$. This is now
comparably easy as we can use the Marsden-Weinstein reduction from
$\mathbb{C}^{n+1}$. We take the constant symplectic structure
$\omega = \frac{\I}{2} g_{\mu\nu} \D z^\mu \wedge \D \cc{z}^\nu$. For
this, the function $g$ is a momentum map for the $\group{U}(1)$
action. Hence $Z$ corresponds to momentum level zero and the
subsequent quotient yields the Marsden-Weinstein reduced phase space
$D_n$. It is now a final but straightforward check that the symplectic
structure obtained this way is compatible with the complex structure
so that we finally end up with a Kähler manifold $D_n$. Note that also
the group $\group{SU}(1, n)$ acts in a Hamiltonian fashion on $D_n$.

We can now construct the star product for the disc as follows: we
start with the Wick star product on $\mathbb{C}^{n+1}$ with respect to
the pseudo Kähler structure defined by $g$ explicitly given by
\begin{equation}
    \label{eq:WickStarProduct}
    f \starwick g
    =
    \sum_{r=0}^\infty \frac{(2\hbar)^r}{r!}
    g^{\mu_1\nu_1} \cdots g^{\mu_r\nu_r}
    \frac{\partial^r f}
    {\partial z^{\mu_1} \cdots \partial z^{\mu_r}}
    \frac{\partial^r g}
    {\partial \cc{z}^{\nu_1} \cdots \partial \cc{z}^{\nu_r}},
\end{equation}
where $g^{\mu\nu} = g_{\mu\nu}$ and summation of the
$\mu_1, \ldots, \nu_r$ is understood as usual. This is a particular
case of the star product for constant Poisson structures on a vector
space as discussed in Section~\ref{sec:WeylStarProduct}. In
particular, we have trivial convergence on polynomials.

To construct the star product on the disc $D_n$ one uses the
$\group{U}(1, n)$-invariance of $\starwick$ and, in particular, the
resulting $\group{U}(1)$-invariance. Classically, the functions
$\Cinfty(D_n)$ can be obtained as the $\group{U}(1)$-invariant
functions on $Z$ and hence as quotient of the $\group{U}(1)$-invariant
functions on $\mathbb{C}^{n+1}$ modulo those which vanish on $Z$.  The
functions vanishing on $Z$ can now be obtained as the ideal generated
by $g-1$ since $1$ is a regular value of $g$. This gives the classical
restriction map
$\Psi_0\colon \Cinfty(\mathbb{C}^{n+1})^{\group{U}(1)} \longrightarrow
\Cinfty(D_n)$
which turns then out to be a $\group{SU}(1, n)$-equivariant Poisson
homomorphism. In a second step, this classical restriction is deformed
into a quantum restriction inducing the star product on the disc.
Instead of formulating this in general as done in
\cite{bordemann.brischle.emmrich.waldmann:1996a,
  bordemann.brischle.emmrich.waldmann:1996b} we give the explicit
formula for this quantum restriction on $\group{U}(1)$-invariant
polynomials.

For multiindices $P$, $Q \in \mathbb{N}_0^{n+1}$ we write
\begin{equation}
    \label{eq:PolynomialdPQ}
    \basis{d}_{P, Q}
    = z^P \cc{z}^Q
    = (z^0)^{P_0} \cdots (z^n)^{P_n}
    (\cc{z}^0)^{Q_0} \cdots (\cc{z}^n)^{Q_n}
\end{equation}
for the usual basis of monomials. Then $\basis{d}_{P, Q}$ is
$\group{U}(1)$-invariant iff $\abs{P} = \abs{Q}$. We denote their
images in $\Cinfty(D_n)$ by
\begin{equation}
    \label{eq:fPQDef}
    \basis{f}_{P, Q} = \Psi_0(\basis{d}_{P, Q}).
\end{equation}
The difficulty is now that the functions $\basis{f}_{P, Q}$ are no
longer linearly independent. Instead, we get relations between them of
the following form. First we specify for multiindices
$P, Q \in \mathbb{N}_0^n$ without zeroth component the special
functions
\begin{equation}
    \label{eq:frPQ}
    \basis{f}_{r, P, Q}
    =
    \begin{cases}
        \basis{f}_{
          (\abs{Q} - \abs{P}, P_1, \ldots, P_n),
          (0, Q_1,\ldots, Q_n)}
        & \textrm{for } \abs{Q} \ge \abs{P} \\
        \basis{f}_{
          (0, P_1, \ldots, P_n),
          (\abs{P} - \abs{Q}, Q_1,\ldots, Q_n)}
        & \textrm{for } \abs{Q} \ge \abs{P}.
    \end{cases}
\end{equation}
They turn out to be linearly independent and we can write every image
$\basis{f}_{P, Q}$ as linear combination
\begin{equation}
    \label{eq:fPQasfrQP}
    \basis{f}_{P, Q}
    =
    \sum_{
      \substack{
        T \in \mathbb{N}_0^n \\
        \abs{T} \le \min\{P_0, Q_0\}
      }
    }
    \binom{\min\{P_0, Q_0\}}{\abs{T}}
    \frac{\abs{T}!}{T!}
    \basis{f}_{r, P'+T, Q'+T},
\end{equation}
where $P' = (P_1, \ldots, P_n) \in \mathbb{N}_0^n$ and analogously for
$Q$. This shows that we obtain a basis for the span of the images of
$\group{U}(1)$-invariant polynomials.

In order to define the quantum restriction map we have to specify the
possible values of $\hbar$ first. While for the
($\group{U}(1)$-invariant) monomials on $\mathbb{C}^{n+1}$ we have
convergence for all $\hbar \in \mathbb{C}$, the quantum restriction
will turn out to be defined only for the following values of
$\hbar$. We define the admissible values for $\hbar$ to be the set
\begin{equation}
    \label{eq:Hdef}
    H
    =
    \mathbb{C} \setminus
    \left\{
        0, -\tfrac{1}{2m}
        \; \big| \;
        m \in \mathbb{N}
    \right\}.
\end{equation}
In particular, the classical limit $\hbar = 0$ is not part of this
open subset but a boundary point. Moreover, we have excluded the
isolated points at $-\tfrac{1}{2m}$ for $m \in \mathbb{N}$. Then we
can define the quantum restriction $\Psi_\hbar$ on the basis of
$\group{U}(1)$-invariant monomials explicitly by
\begin{equation}
    \label{eq:QuantumRestriction}
    \Psi_\hbar(\basis{d}_{P, Q})
    =
    (2\hbar)^{\abs{P}}
    \bigg(\frac{1}{2\hbar}\bigg)_{\abs{P}}
    \Psi_0(\basis{d}_{P, Q})
    =
    (2\hbar)^{\abs{P}}
    \bigg(\frac{1}{2\hbar}\bigg)_{\abs{P}}
    \basis{f}_{P, Q},
\end{equation}
where $\hbar \in H$ and $(z)_m = z(z+1) \cdots (z+m-1)$ is the
Pochhammer symbol. The relevance of this definition is now the
following:
\begin{lemma}
    \label{lemma:RestrictionKernel}%
    The kernel of the linear map $\Psi_\hbar$ coincides with the
    two-sided $\starwick$-ideal inside the $\group{U}(1)$-invariant
    polynomials generated by $g-1$.
\end{lemma}
The original construction of
\cite{bordemann.brischle.emmrich.waldmann:1996a,
  bordemann.brischle.emmrich.waldmann:1996b} did not use the
polynomial functions to construct the reduction, resulting in a
slightly more complicated definition. However, there the ultimate
formula was derived from a more conceptual point of view instead of
our short-cut by guessing \eqref{eq:QuantumRestriction} and verifying
its properties afterwards. When working only with polynomials as we do
here, the verification of the lemma is a fairly simple computation.

Thanks to this observation we can define a product on the span of the
image of $\Psi_0$, i.e. for the span of the functions
$\basis{f}_{P, Q}$ simply by pushing forward $\starwick$. The algebra
obtained this way is the quotient of the $\group{U}(1)$-invariant
polynomials modulo the $\starwick$-ideal generated by
$g-1$. Explicitly, this results in the formula
\begin{equation}
    \label{eq:ProductForfPQs}
    \basis{f}_{P, Q} \starred \basis{f}_{R, S}
    =
    \sum_{T = 0}^{\min\{P, S\}}
    (-1)^{T_0}
    \frac{
      (\frac{1}{2\hbar})_{\abs{P+S-T}} T!
    }{
      (\frac{1}{2\hbar})_{\abs{P}}
      (\frac{1}{2\hbar})_{\abs{S}}
    }
    \binom{P}{T}\binom{S}{T}
    \basis{f}_{P+R-T, Q+S-T}
\end{equation}
for all $P, Q, R, S \in \mathbb{N}_0^{1+n}$ with $\abs{P} = \abs{Q}$
and $\abs{R} = \abs{S}$. Note however, that the functions
$\basis{f}_{P, Q}$ do not yet form a basis, one has to expand them
according to \eqref{eq:fPQasfrQP} using the basis of the functions
$\basis{f}_{r, P, Q}$ with $P, Q \in \mathbb{N}_0^n$ instead.

Now we can use the continuity estimates for $\starwick$ from
Section~\ref{sec:WeylStarProduct}. Since we have here a
finite-dimensional situation it is convenient to describe the
seminorms $\halbnorm{p}_R$ defining the $\Sym_R$-topology more
explicitly. We focus on the limiting case $R = \frac{1}{2}$ directly.
Let $\rho > 0$. Then we define for a polynomial $a$ on
$\mathbb{C}^{n+1}$ the norm
\begin{equation}
    \label{eq:NormRho}
    \norm[\bigg]{
      \sum_{P, Q \in \mathbb{N}_0^{n+1}}
      a_{P, Q} \basis{d}_{P, Q}
    }_{\mathbb{C}^{n+1}, \rho}
    =
    \sum_{P, Q \in \mathbb{N}_0^{n+1}}
    \abs{a_{P, Q}} \rho^{\abs{P + Q}} \sqrt{\abs{P + Q}!}.
\end{equation}
Note that for a polynomial the series is in fact a finite sum.
\begin{lemma}
    \label{lemma:SameTopology}%
    The collection $\norm{\argument}_{\mathbb{C}^{n+1}, \rho}$ of
    norms for $\rho > 0$ defines the $\Sym_R$-topology for the
    polynomials on $\mathbb{C}^{n+1}$ for $R = \frac{1}{2}$.
\end{lemma}

For the functions on the disc we consider the basis $\basis{f}_{r, P,
  Q}$ and specify norms on the span of these functions, i.e. on the
quotient algebra. We set
\begin{equation}
    \label{eq:NormForDisc}
    \norm[\bigg]{
      \sum_{P, Q \in \mathbb{N}_0^n}
      a_{P, Q} \basis{f}_{r, P, Q}
    }_{D_n, \rho}
    =
    \sum_{P, Q \in \mathbb{N}_0^n}
    \abs{a_{P, Q}} \rho^{\abs{P + Q}}.
\end{equation}
for $\rho > 0$. By a direct estimate one can then show that the
resulting locally convex topology for the span of the
$\basis{f}_{r, P, Q}$ coincides with the quotient topology induced by
$\Psi_\hbar$:
\begin{lemma}
    \label{lemma:QuotientTopology}%
    The locally convex topology defined by the norms
    $\norm{\argument}_{D_n, \rho}$ for $\rho > 0$ coincides with the
    quotient topology induced by $\Psi_\hbar$.
\end{lemma}
In particular, the star product $\starred$ becomes continuous with
respect to this topology. One can find explicit norm estimates for
$\starred$ but the next theorem already follows already from the
previous two lemmas and our investigations concerning $\starwick$ in
Section~\ref{sec:WeylStarProduct}.
\begin{theorem}
    \label{theorem:DiscContinuous}%
    Let $\hbar \in H$. Then the star product $\starred$ for the span
    of the functions $\basis{f}_{r, P, Q}$ is continuous with respect
    to the quotient topology explicitly described by the norms
    \eqref{eq:NormForDisc}.
\end{theorem}
In \cite{beiser.waldmann:2014a} a slightly finer topology was used,
the above one seems to be more appropriate as it will allow for a
larger completion. In fact, the completion in
\cite{beiser.waldmann:2014a} was identified with a certain nuclear
Köthe space. Nevertheless, the interpretation as functions was still
not satisfactory. The above topology will now allow to determine the
functions in the completion very explicitly and geometrically.

To describe the completion we first note that all the functions
$\basis{d}_{P, Q}$ as well as $\basis{f}_{P, Q}$ are real-analytic of
a very particular form. Being real-analytic means that they can be
extended to an open neighbourhood of the diagonal in the Cartesian
product of the underlying Kähler manifold with itself. Now the
polynomials as well as the $\basis{f}_{r, P, Q}$ can be extended
holomorphically to a much larger complex manifold than an open
neighbourhood of the diagonal. To describe these new complex
manifolds we have to double all previous ones as follows:

First we consider $\mathbb{C}^{n+1} \times \mathbb{C}^{n+1}$ with the
$\group{U}(1, n)$-action defined by $U \acts (x, y) = (Ux, \cc{U}y)$.
This way, we have the $\group{U}(1, n)$-equivariant diagonal map
$\Delta(p) = (p, \cc{p})$. Let
$\tau\colon \mathbb{C}^{n+1} \times \mathbb{C}^{n+1} \longrightarrow
\mathbb{C}^{n+1} \times \mathbb{C}^{n+1}$
be the anti-holomorphic flip diffeomorphism defined by
$\tau(x, y) = (\cc{y}, \cc{x})$. Then the image of $\Delta$ coincides
with those points $(x, y)$ with $\tau(x, y) = (x, y)$, i.e. the fixed
points of $\tau$. We can now extend the function $g$ to a holomorphic
function
\begin{equation}
    \label{eq:HolomorphicHatg}
    \hat{g}
    =
    g_{\mu\nu} x^\mu y^\nu
    \in
    \Holomorph(\mathbb{C}^{n+1} \times \mathbb{C}^{n+1}).
\end{equation}
Then $\hat{g} \circ \Delta = g$ and thus $\hat{g}$ is the unique
holomorphic extension of the real-analytic function $g$.
Using $\hat{g}$ we can define the complex submanifold
\begin{equation}
    \label{eq:hatZ}
    \hat{Z}
    =
    \hat{g}^{-1}(\{1\})
    \subseteq
    \mathbb{C}^{n+1} \times \mathbb{C}^{n+1},
\end{equation}
which is indeed a complex submanifold as $1$ is a regular value of
$\hat{g}$ which is invariant under $\tau$. Clearly, $\group{U}(1, n)$
acts on $\hat{Z}$. Moreover, the diagonal $\Delta$ maps $Z$ into
$\hat{Z}$ and the image consists of the fixed points in $\hat{Z}$
under $\tau$.

In a next step we extend the $\group{U}(1)$ action on $Z$ to an action
of the multiplicative Lie group
$\mathbb{C}_* = \mathbb{C} \setminus \{0\}$ on $\hat{Z}$ as
follows. For $z \in \mathbb{C}_*$ we define
$z \acts (x, y) = (zx, \frac{1}{z}y)$ resulting in a holomorphic
action. The action commutes with the $\group{U}(1, n)$ action and
provides a free and proper action on $\hat{Z}$. Thus the quotient
\begin{equation}
    \label{eq:hatD}
    \hat{D}_n = \hat{Z} \big/ \mathbb{C}_*
\end{equation}
becomes a complex manifold which inherits the $\group{U}(1, n)$
action. Of course, only the subgroup $\group{SU}(1, n)$ acts
non-trivially. The anti-holomorphic involution $\tau$ also descends to
an anti-holomorphic involution $\tau$ on $\hat{D}_n$.  The original
disc $D_n$ can now be included into $\hat{D}_n$ using the diagonal
$\Delta$. Indeed, the diagonal inclusion of $Z$ into $\hat{Z}$
descends to a smooth inclusion
$\Delta_D\colon D_n \longrightarrow \hat{D}_n$ as a real submanifold,
compatible with the $\group{U}(1, n)$ actions. Note that $\hat{D}_n$
can also be seen as a open subset of
$\mathbb{CP}^n \times \mathbb{CP}^n$.

It will be this doubled version of the disc encoding the
completion. First we note the following observation:
\begin{lemma}
    \label{lemma:RestrictHolToDiag}%
    A holomorphic function $\hat{a} \in \Holomorph(\hat{D}_n)$ is
    uniquely determined by its restriction
    $\Delta_D^*\hat{a} \in \Cinfty(D_n)$.
\end{lemma}
This allows to consider those (necessarily real-analytic) functions on
$D_n$ which are restrictions of holomorphic functions on $\hat{D}_n$,
i.e. those which allow for a (necessarily unique) holomorphic
extension to $\hat{D}_n$. We define
\begin{equation}
    \label{eq:AlgebraOnDn}
    \algebra{A}(D_n)
    =
    \left\{
        a \in \Cinfty(D_n)
        \; \Big| \;
        a = \Delta_D^* \hat{a}
        \textrm{ for some }
        \hat{a} \in \Holomorph(\hat{D}_n)
    \right\}.
\end{equation}
Then
$\Delta_D^*\colon \Holomorph(\hat{D}_n) \longrightarrow
\algebra{A}(D_n)$
becomes an algebra isomorphism with respect to the commutative
product. This allows to pull the standard locally convex topology of
holomorphic functions back to $\algebra{A}(D_n)$. Explicitly, we
define seminorms of functions in $\algebra{A}(D_n)$ by
\begin{equation}
    \label{eq:SupKNorms}
    \norm{a}_{D_n, K}
    =
    \sup_{u \in K} \abs{\hat{a}(u)}
\end{equation}
for a compact subset $K \subseteq \hat{D}_n$, where $\hat{a}$ is the
unique holomorphic function with $\Delta_D^*\hat{a} = a$ as
before. This way, $\Delta_D^*$ becomes an isomorphism of Fréchet
spaces. The topology is the Fréchet topology of locally uniform
convergence on $\hat{D}_n$.

Since the polynomials $\basis{d}_{P, Q}$ can be extended to
holomorphic polynomials on $\mathbb{C}^{n+1} \times \mathbb{C}^{n+1}$
it follows fairly easy that the functions $\basis{f}_{P, Q}$ can be
extended holomorphically to $\hat{D}_n$, i.e. we have
\begin{equation}
    \label{eq:fPQExtension}
    \basis{f}_{P, Q} = \Delta_D^* \hat{\basis{f}}_{P, Q}
    \in \algebra{A}(D_n)
\end{equation}
with some unique $\hat{\basis{f}}_{P, Q} \in \Holomorph(\hat{D}_n)$.
In particular, this applies to the basis $\basis{f}_{r, P, Q}$ itself.
It requires now an adaption of the Cauchy integral formula to see that
the corresponding functions
$\hat{\basis{f}}_{r, P, Q} \in \Holomorph(\hat{D}_n)$ play the role of
monomials: more precisely, they form an absolute Schauder basis with
coefficient functionals given by iterated Cauchy integrals. This means
that every holomorphic function $\hat{a} \in \Holomorph(\hat{D}_n)$
has a unique absolutely convergent expansion
\begin{equation}
    \label{eq:hataExpansion}
    \hat{a}
    =
    \sum_{P, Q \in \mathbb{N}_0^n}
    \hat{a}_{P, Q} \hat{\basis{f}}_{r, P, Q}
    \quad
    \textrm{with}
    \quad
    a_{P, Q}
    =
    \frac{1}{(-4\pi^2)^n}
    \oint \cdots \oint
    \hat{a}
    \frac{(1 - uv)^{\max\{\abs{P}, \abs{Q}\}}}
    {v^{P+1}u^{Q+1}} \D^nu \D^nv,
\end{equation}
with the canonical holomorphic coordinates $(u, v)$ of $\hat{D}_n$
inherited from the standard chart of
$\mathbb{CP}^n \times \mathbb{CP}^n$ where we just divide by the
zeroth component as usual. The contour integrals are around $(0, 0)$
in this chart. Transferring this back to $\algebra{A}(D_n)$ gives the
same expansion for every $a \in \algebra{A}(D_n)$ with the very same
coefficients.

The last step consists now in recognizing that the topology on the
polynomial functions for which $\starred$ was shown to be continuous
simply \emph{coincides} with the canonical topology of the holomorphic
functions. This ultimately yields the following description of the
completion:
\begin{theorem}
    \label{theorem:CompletionDisc}%
    Let $\hbar \in H$. The completion of the span of the functions
    $\basis{f}_{r, P, Q}$ with respect to the quotient topology,
    i.e. the seminorm system given by \eqref{eq:NormForDisc},
    coincides with $\algebra{A}(D_n)$ and the resulting Fréchet
    topology is the inherited one from $\Holomorph(\hat{D}_n)$. Thus
    $\starred$ extends to $\algebra{A}(D_n)$.
\end{theorem}
\begin{remark}[Further properties]
    \label{remark:Properties}%
    Let $\hbar \in H$.
    \begin{remarklist}
    \item \label{item:Expansion} First we note that the expansion
        \eqref{eq:hataExpansion} converges in the
        $\algebra{A}(D_n)$-topology. Since the product is continuous,
        the product of two elements $a, b \in \algebra{A}(D_n)$ in the
        completion can be computed by the absolutely convergent series
        \begin{equation}
            \label{eq:astarredb}
            a \starred b
            =
            \sum_{P, Q, R, S \in \mathbb{N}_0^n}
            a_{P, Q} b_{R, S}
            \basis{f}_{r, P, Q} \starred \basis{f}_{r, R, S}
        \end{equation}
        with the remaining product of the basis functions given by
        \eqref{eq:ProductForfPQs}. This gives a very explicit way to
        compute the star product on the disc.
    \item \label{item:Involution} Since the original Wick star product
        had the pointwise complex conjugation as $^*$-involution if
        $\hbar = \cc{\hbar}$ is real, it follows easily that also
        $\starred$ has the pointwise complex conjugation as
        (continuous) $^*$-involution. This turns
        $(\algebra{A}(D_n), \starred)$ into a Fréchet $^*$-algebra.
    \item \label{item:HolomorphicDef} For any two elements $a, b \in
        \algebra{A}(D_n)$ the map
        \begin{equation}
            \label{eq:hbarDependence}
            H \ni \hbar
            \; \mapsto \;
            a \star_{{\scriptscriptstyle\mathrm{red}}, \hbar} b
            \in \algebra{A}(D_n)
        \end{equation}
        is a holomorphic map with respect to the Fréchet topology of
        $\algebra{A}(D_n)$. Here we write explicitly
        $\star_{{\scriptscriptstyle\mathrm{red}}, \hbar}$ to emphasize
        the $\hbar$-dependence of the reduced star product. One can
        construct explicit elements such that this statement can not
        be improved: there are $a$ and $b$ such that $a \starred b$
        has poles at \emph{every} critical value
        $\hbar = - \frac{1}{2m}$. Note that for the original
        polynomials only finitely many poles occurred in a given
        product, depending on the maximal degree.
    \end{remarklist}
\end{remark}

It is the last remark which makes the analysis of the semi-classical
limit more complicated as the nice cases discussed in
Section~\ref{sec:WeylStarProduct} and
Section~\ref{sec:GuttStarProduct}. Now the classical limit $\hbar = 0$
is not an interior point of the domain where the product depends
holomorphically on $\hbar$. It is a boundary point with poles
accumulating on the negative axis. Thus the limit
$\hbar \longrightarrow 0$ has to be taken with care from the right,
i.e. we only can expect a reasonable limit for
$\hbar \longrightarrow 0^+$. The following theorem shows that this can
be done.
\begin{theorem}
    \label{theorem:SemiclassicalDisc}%
    Let $a, b \in \algebra{A}(D_n)$. Then one has
    \begin{equation}
        \label{eq:SemiClassicalLimit}
        \lim_{\hbar \to 0^+}
        a \star_{{\scriptscriptstyle\mathrm{red}}, \hbar} b
        =
        ab
        \quad
        \textrm{and}
        \quad
        \lim_{\hbar \to 0^+}
        \frac{\I}{\hbar}
        (a \star_{{\scriptscriptstyle\mathrm{red}}, \hbar} b
        -
        b \star_{{\scriptscriptstyle\mathrm{red}}, \hbar} a)
        =
        \{a, b\}.
    \end{equation}
\end{theorem}
The above statement becomes trivial for finite linear combinations in
the basis functions $\basis{f}_{r, P, Q}$ as there the reduced star
product is still holomorphic around $\hbar = 0$ and there are only
finitely many poles on the negative axis. However, the above example
shows that for general $a$ and $b$ the argument is more subtle. In
fact, its proof is fairly technical and requires a considerable
effort.

We conclude this section with a few remarks on further properties of
the star product $\starred$. In
\cite{kraus.roth.schoetz.waldmann:2018a:pre} it was shown that even
though the topology is not locally multiplicatively convex for the
star product, many interesting transcendental functions can be
defined. In particular, the algebra $\algebra{A}(D_n)$ contains many
bounded functions even though the span of the Schauder basis only
contains constant bounded functions. The reason is simply that we have
real functions $a \in \algebra{A}(D_n)$ which can then be
exponentiated to yield bounded functions
$\Delta_D^* \exp(\I \hat{a}) \in \algebra{A}(D_n)$. This gives some
hope that we also have \emph{periodic} functions with respect to
Fuchsian groups. In this case, one would obtain immediately a
construction of convergent star products on Riemann surfaces of higher
genus. Ultimately, one can then compare this approach with the one of
Bieliavsky \cite{bieliavsky:2017a}.

For $\hbar > 0$ it was shown that every classically positive linear
functional of $\Cinfty(D_n)$, i.e. a positive Radon measure on $D_n$,
is also positive for the $^*$-algebra $(\algebra{A}(D_n), \starred)$.
In particular, all evaluation functionals at the point of the disc
$D_n$ are positive functionals. This shows that the resulting star
product algebra has sufficiently many positive linear functionals to
separate elements. In particular, it has faithful $^*$-representations
on pre-Hilbert spaces, see \cite{schmuedgen:1990a} for more details on
the $^*$-representation theory of unbounded operator algebras and
\cite{waldmann:2005b, bursztyn.waldmann:2005a, waldmann:2019a:script,
  bursztyn.waldmann:2005b} for some algebraic background.

The components of the momentum map for the $\group{SU}(1, n)$ are
contained in $\algebra{A}(D_n)$: they are linear combinations of the
Schauder basis $\basis{f}_{r, P, Q}$ for small $P$ and $Q$. However,
their $\starred$-exponentials are not defined as elements of
$\algebra{A}(D_n)$. Nevertheless, suppose that one has a
$^*$-representation $\pi$ of $\algebra{A}(D_n)$ for $ \hbar > 0$ on
some pre-Hilbert space $\mathcal{H}$ such that
$a \mapsto \abs{\SP{\phi, \pi(a)\psi}}$ is continuous for all
$\phi, \psi \in \mathcal{H}$. Here every GNS representation for a
continuous positive linear functional gives an example. Then one can
show that the representing operators for the components of the
momentum map are essentially self-adjoint and integrate to a strongly
continuous unitary representation of $\group{SU}(1, n)$ on the Hilbert
space completion of $\mathcal{H}$. This will allow to use the star
product algebra to explicitly construct well-behaved representations
of $\group{SU}(1, n)$.

%
% Open Questions and Outlook
%

\section{Open Questions and Outlook}
\label{sec:OpenQuestionsOutlook}

Let us conclude this review with some open questions and hints on a
further development for the theory of convergence of star products.
\begin{enumerate}
\item The three classes of examples indicate that the proposed way to
    convergence in Section~\ref{sec:QuestConvergence} is at least not
    completely hopeless. Of course, the examples are the most simple
    ones but already here we see a rich and nontrivial analytic
    structure when discussing convergence. Thus it is reasonable to
    stick with the proposal and investigate further examples.
\item The current situation of convergence of star products based on a
    detailed analysis of the formal power series should of course be
    compared to the previous approaches based on integral
    formulas. Here the situation is not completely easy as the
    functional spaces are typically rather complementary. On the one
    hand one poses conditions on the growth of the Taylor coefficients
    but allows fairly unbounded behaviour at infinity, on the other
    side one can abandon real-analytic functions and allow for more
    general smooth functions as long as the growth at infinity is
    either bounded or at least very moderate. Nevertheless, a
    comparison should be possible. In particular, when passing to
    suitable $^*$-representations then a comparison by means of
    spectral measures is one option.
\item The formal power series approach is perhaps the only possibility
    to include infinite-dimensional phase spaces. The applications to
    (quantum) field theory are by far not yet explored. Here already
    the constant Poisson structures provide interesting scenarios, see
    e.g. \cite[Sect.~6]{waldmann:2014a}.
\item The above approaches are always based on the assumption to have
    a continuous product. As long as one is heading for Fréchet
    topologies this is not a big restriction. However, in various
    field theoretic models one might expect that already the Poisson
    bracket has less regularity: it can be only separately continuous
    or even only sequentially separately continuous. It remains an
    open question to develop appropriate tools to discuss the
    convergence for formal star products quantizing such Poisson
    structures.
\item At the moment it is perhaps too early to propose general
    definitions for convergence schemes of formal star products. In
    particular, the example of the disc shows that the dependence on
    $\hbar$ might be more delicate than a first guess suggest. It
    seems that more examples have to be investigated. Here the Wick
    type star products might turn out to be very handy as they provide
    the better positivity properties. Having a real-analytic context
    for the manifold directions will also help to understand the
    convergence of the formal power series in $\hbar$. Thus
    non-compact Hermitian symmetric spaces will be natural candidates
    for further investigations.
\item The convergence properties of the above examples can be studied
    further and will lead to applications beyond deformation
    quantization. In particular, the Gutt star product suggest
    immediately next steps in direction of representation theory as
    the self-adjointness of Lie algebra representations can now be
    studied with the help of the completed algebra. E.g. in GNS
    representations one can use the completed and therefore rich
    algebra to obtain analytic vectors and hence self-adjointness. But
    also within differential geometry there will be applications when
    thinking of the above Fréchet algebras as algebras of particular
    pseudo-differential operators. With the help of
    $^*$-representations, one has the possibilities for spectral
    analysis in these algebras.
\item Finally, in the theory of locally convex algebras much effort is
    spend on the cases of locally multiplicatively convex
    algebras. Here one can use Banach algebra techniques and obtains a
    closely related extension compared to the Banach algebra
    case. However, all the above examples indicate that there is a
    non-trivial and interesting world of locally convex algebras far
    beyond the locally multiplicatively convex case. Here one needs to
    develop new techniques also from a more conceptual point of
    view. The above examples can then be seen as a guideline to
    formulate such an extension: they should provide ideal testing
    grounds as they show already a fairly complicated nature but are
    still manageable and concrete enough to test general ideas.
\end{enumerate}

%
% The bibliographies, uncomment if you need to cite someone
% Make it smaller than other text
%

{
  \footnotesize
%  \bibliographystyle{chairx}
%  \bibliography{dqarticle,dqbook,dqprocentry,dqproceeding,dqthesis,preprints,misc,notes,script}

}

%
% this is the end of constar
%

\end{document}

%%% Local Variables:
%%% mode: latex
%%% TeX-master: t
%%% End: